# Roughness in BO/BH/Z-ALGEBRA

## Faraj.A.Abdunabi[1], Ahmed shletiet[2]

[1] (Department of Mathematics, Faculty of Science, AJdabyia University, Libya.)
Email: faraj.a.abdunabi@uoa.edu.ly
[2] (Department of Mathematics, Faculty of Science, AJdabyia University, Libya.)
Email:Ahmed.shletiet@uoa.edu.ly

**الملخص**

الهدف الرئيسي من هذه الورقة هو تقديم مفاهيم الجبر الخشن -Z / BH / BO كتوسيع لمفهوم -Z / BH / BO الجبر على التوالي. الهدف الأخر هو النظر في الرواسم (القوي) ذو القيمة المحددة في هذه الهياكل الجبرية. يتم التحقيق في مفهوم التشكل Z / BH / BO (القوي) ذو القيمة المحددة في جبر Z / BH / BO بعدة خصائص. باستخدام مفهوم مساحة التقريب المعمم والمثالية لجبر BO / BH / Z ، ندرس مفهوم كنوعًا آخرن من التقريبات من أعلي ومن أسفل المعممة بناءً على النموذج المثالي. بالإضافة إلى ذلك ، تمت دراسة بعض الخصائص.

## *Abstract*

The main goal of this paper, present the concepts of rough BO/BH/Z- Algebra as extended of the concept of BO/BH/Z-algebra respectively. The other goal is to consider the (strong) set-valued mapping in these algebraic structures. The concept of a (strong) set-valued BO/BH/Z-morphism in BO/BH/Z algebras is investigated with several properties. Using the concept of generalized approximation space and ideal of BO/BH/Z-algebra, we consider another type of generalized lower and upper approximations based on the ideal. In addition, some properties are studied.

## 1. Introduction

TheRough set theory was present by Pawlak [1] in 1982. It is a good tool for modeling and processing incomplete information in the information system. The concepts of rough set theory build of lower and upper approximations. J. Neggers and H. S. Kim [2] introduce the concept of B-algebras. In[3], Young. B J. and el. consider the fuzzification of (normal) B-subalgebras in B-algebras. In[4] Chang Bum Kim and Hee Sik Kim introduce the notion of a BO-algebra. Y. B. Junand et[5] introduced the concept of a BH-algebra. The Z-algebra present by M. Chandramouleeswaran And Et.in[6].The main purpose of this paper is to introduce rough BO/BH/Z-algebra as extended of the concept of BO-algebra (BH-algebra ) respectively Moreover, we introduce some properties of approximations and these algebraic structures.

## 2. Preliminaries

We start by giving some definitions and results about rough sets.

Suppose that R is an equivalence relation on a universe set (nonempty finite set) U. The pair (U,R) is denoted to the approximation space. The notation U/R is denoted as the family of all equivalent classes[a]R. The empty set is $\varnothing$, the elements of U/R are called elementary sets, and Ac is a complementation of A For any $A \subseteq U$.

**Definition 2.1**: Let (U, R) be an approximation space. Define the upper approximation of A is $\overline{RA} = \{a \in U : [a]_R \cap A \neq \emptyset\}$ and the lower approximation of A is $\underline{RA} = \{a \in U : [a]_R \subseteq A\}$ the boundary is $BA_R = \overline{RA} - \underline{RA}$. If $BA_R = \emptyset$, then A is the exact (crisp) set, and if $BA_R \neq \emptyset$, X is a rough set ( inexact).

**Preposition 2-1**: Suppose that (U, R) is an approximation space. Let A,B$\subseteq$U, then:

1) $\underline{RA} \subseteq A \subseteq \overline{RA}$,
2) $\underline{R\emptyset} = \overline{R\emptyset}, \underline{RU} = \overline{RU}$,
3) $\underline{R(A \cup B)} \supseteq \underline{R(A)} \cup \underline{R(B)}$,
4) $\underline{R(A \cap B)} = \underline{R(A)} \cap \underline{R(B)}$,
5) $\overline{R(A \cup B)} = \overline{R(A)} \cup \overline{R(B)}$.
6) $\overline{R(A \cap B)} \subseteq \overline{R(A)} \cap \overline{R(B)}$.
7) $\overline{RA^C} = \left(\underline{RA}\right)^c$.
8) $\underline{RA^C} = (\overline{RA})^c$.
9) $\underline{R(\underline{RA})} = \overline{R(\underline{RA})} = \underline{RA}$.
10) $(R(\overline{RA}) = \underline{R(\overline{RA})} = \overline{RA}$.
11) $\overline{RA} \ \overline{RB} = \overline{RAB}$.
12) $\underline{RA} \ \underline{RB} \subseteq \underline{RAB}$.

The concept of BO/ BH/Z-algebra with examples are discussed in this portion.

**Definition 2.2:** Let X be a non-empty set with binary process $*$, $0 \in X$ is B-algebra if $\forall$ x, y, z$\in$ X sitsifies:

C1: $x * x = 0$.

C2: $x * 0 = x$.

C3 :$(x*y)* z = x*(z*(0*y))$.

where 0 is called zero element.

**Remark 2.1.** The element $e \in X$ is called right-identity if x*e =x and left identity if e*x=x for every x∈X and x≠e. e is called the identity if x*e =x and e*x=x for every x∈X. Then (X,*) is called B-algebra containing identity.

**Example 2.1.** Suppose that X = {0,1,2 e}. Define the binary operation on X as shown in the following table 1

| * | 0 | 1 | 2 | e |
|---|---|---|---|---|
| 0 | 0 | 1 | 2 | e |
| 1 | 1 | 0 | e | 2 |
| 2 | 2 | e | 0 | 1 |
| e | e | 2 | 1 | 0 |

Table 1

Table 1 shows that the (X, *) is B-algebra with the identity element.

**Definition 2.3.** Let X be a non-empty set with binary process *, 0∈X is BH-algebra if ∀ x, y, z∈ X sitsifies:

(C1), (C2), and

(C4) $For\ any\ x, y\ \in\ X, x * y = y * x = 0 \Rightarrow x = y$.

**Definition 2.4** Let X be a non-empty set with binary process *, 0∈X is BO-algebra if ∀ x, y, z∈ X sitsifies:

(C1), (C2) and

C5: $x * (y * z) = (x * y) * (0 * z)\ for\ any\ x, y, z\ \in\ X$.

**Example 2.2:** Suppose that X={ 0,1,2,3,4} and the following table 2 of * :

| * | 0 | 1 | 2 | 3 | 4 |
|---|---|---|---|---|---|
| 0 | 0 | 2 | 1 | 4 | 3 |
| 1 | 1 | 0 | 3 | 2 | 4 |
| 2 | 2 | 4 | 0 | 3 | 1 |
| 3 | 3 | 1 | 4 | 0 | 2 |
| 4 | 4 | 3 | 2 | 1 | 0 |

Table 2

Table 2 shows that the (X,*,0) is BO-algebra.

**Example 2.3.** Suppose that X={ 0,1,2 ,3} and the following table 3 of *

| * | 0 | 1 | 2 | 3 |
|---|---|---|---|---|
| 0 | 0 | 1 | 0 | 0 |
| 1 | 1 | 0 | 0 | 0 |
| 2 | 2 | 2 | 0 | 3 |
| 3 | 3 | 3 | 3 | 0 |

Table 3

Table 3 shows that the (X,*,0) is a BH-algebra.

**Definition 2.4**. Suppose (I≠φ) ⊆ BH/Z)-algebra. I is called a BH/Z-ideal of X respectivtly if it satisfies the following conditions:

(1) $0 \in I$,

(2) $(x * y) \in I, y \in I \Rightarrow x \in I, \forall x, y, z \in X$.

(3) $(x * y) * z \in I, y \in I \Rightarrow x * z \in I, \forall x, y, z \in X$ , then I called strong Ideal of X.

**Definition 2.5[6]**. Let X be a non-empty set with binary process *, 0∈X is Z-algebra if ∀ x, y, z∈ X sitsifies(C1-C2) and

C6:$x * x = x$

C7: $x * y = y * x$, when $x \neq 0$ and $y \neq 0, \forall x, y \in X$.

**Example 2.4**. Suppose that X={ 0,1,2,3 } and the following table 4 of *

| * | 0 | 1 | 2 | 3 |
|---|---|---|---|---|
| 0 | 0 | 1 | 2 | 3 |
| 1 | 1 | 0 | 0 | 1 |
| 2 | 0 | 0 | 2 | 2 |
| 3 | 0 | 1 | 2 | 3 |

Table 4

Table 4 shows that the (X,*,0) is a Z-algebra. If I = {0, 1, 2}, then it is is a Z-ideal of X.

**3. Main Result**

**Definition 3.1**: Suppose that ~ be an equivalence relation on a set X=(X,*,0). If x ∈ X, defined [x]~ the ~class of x s follows: [x]~ = {y ∈ X | (x, y) ∈ ~}.The equivalence relation ~ on X is called a congruence relation if

$(\forall x, y, z \in X) ((x, y) \in \sim \Rightarrow (x * y, y * z) \in \sim, (z * x, z * y) \in \sim)$ .

**Definition 3.2.** Suppose that A and B two non-empty subsets of X, we denote $AB = A * B = \{a * b \,|\, a \in A \text{ and } b \in B\}$. Let $\sim$ be an equivalence relation on X. Then $(\forall x, y \in X)([x]\sim[y]\sim \subseteq [x * y]\sim)$.

If $Y \in P(X)$, we define the upper approximation of Y by $+[Y]\sim = \{x \in X \,|\, [x]\sim \subseteq Y\}$ and the lower approximation of Y is $-[Y]\sim = \{x \in X \,|\, [x]\sim \cap Y \neq \emptyset\}$. The pair $(X, \sim)$ is called an approximation space.

Note that, $+[Y]\sim$ and $-[Y]\sim$ are subsets of X.

If $Y \subseteq X$, then Y is said to be definable if $+[Y]\sim = -[Y]\sim$ and rough otherwise.

Suppose that I be a BO/BH-ideal of X. Define a relation $\sim$ on X by $(x, y) \in \sim$ if and only if $x* y \in$ I and $y*x \in$ I.

**Definition 3.3.** Suppose that $(X, \sim)$ is an approximation space, a pair $(I1, I2) \in P(X) \times P(X)$ is called a rough set in $(X,\sim)$ if and only if $(I1, I2) = Apr(X)$ for some $X \in P(X)$.

**Example 3.1**: consider example 2.2. Let $Y = \{0, 1\}$ be a BO-ideal of X. Suppose that $\sim$ is an equivalence relation on X related to Y.

So, Y0 = Y1 = Y, Y2 = {2}, Y3 = {3}, and Y4 = {4}. Hence, -[Y, {0,1}] = {0, 1} , -[Y, {0,2}] = {2}, -[Y, {0,3}]= {3}, and -[Y, {0,1,2,3}] = {0, 1, 2,3}. However, + [Y, {0,1}] = {0, 1}.+ [Y, {0}] = {0, 1}, + [Y, {2}] = {0, 2} , + [Y, {1,2,3}] = {0,1, 2, 3} , + [Y, {0,2,3}] = {0, 1, 2, 3},+ [Y, {1,2,3,4}] = {0, 1, 2,3, 4}.

Here, there exists a non-BO-ideal Y of X such that their lower and upper approximation are BO-ideals of X.

**Proposition 3.1.** Let X be a Bo(BH)-algebra and A, B two subsets of X. Let $\sim$ be an equlivence relations on X. Then the following hold:

1) $-[A]\sim \subseteq A \subseteq +[A]\sim,$
2) $+[A \cup B]\sim = +[A]\sim \cup +[B]\sim,$
3) $-[A \cap B]\sim = -[A]\sim \cap -[B]\sim,$
4) $If\ A \subseteq B, then -[A] \sim \subseteq -[B] \sim and +[A] \sim \subseteq +[B] \sim,$
5) $-[A]\sim \cup -[B]\sim \subseteq -[A \cup B]\sim,$
6) $+[A \cap B] \sim \subseteq +[A] \sim \cap +[B] \sim.$

Proof. Straightforward.

Let X be a BH-algebra and let $\emptyset \neq$ A, B $\subseteq$ X. Define A $* $ B := {a $*$ b|a $\in$ A, b $\in$ B}.

**Proposition 3.2.[7]**. Suppose that X is BH-algebra. Let ~ be a congruence relation on X. Suppose that A, B are two non-subsets of X. Then

1) $+[A]\sim * +[B]\sim \subseteq +[A*B]\sim$.
2) $If -[A*B] \neq \phi, then -[A]\sim * -[B]\sim \subseteq -[A*B]\sim$.

Proof.

Assum that x ∈+[A]~*+[B] ~. Then x = a∗b for some a ∈ +[A]~ and b ∈+[B].

Then, we have y, z ∈ X such that y∈[a]~∩A and z∈[b]~∩B. Hence y∈[a]~, z∈ [z] ~, y∈A and z∈B. Since ~ is a congruence relation on X, y∗z ∈ [a]~∗[b]~= [a∗b]ρ. Since y∗z ∈ A ∗ B, we have x = a ∗ b ∈+[A∗B].

Suppose that x∈−[A]~∗−[B]~. Then x = a∗b for some a ∈ −[A]~ and b∈−[B]~. Thus we have [a]~⊆ A and [b]~⊆ B. [a∗ b]~ = [a]~∗[b]~ ⊆ A ∗ B because ~ is a congruence relation on X. Then, x =a∗b ∈ −[A∗B] ~.

**Proposition 3.3**. Suppose that X is BO/Z-algebra. Let ~ be a congruence relation on X. Suppose that A, B are two non-subsets of X. Then

1) $+[A]\sim * +[B]\sim \subseteq +[A*B]\sim$.
2) $If -[A*B] \neq \phi, then -[A]\sim * -[B]\sim \subseteq -[A*B]\sim$.

Proof the same strategy in Proposition 3.2.

**Definition 3.4**. Let X and Y be non-empty universes and consider the mapping F : X → P(Y). we say F is a set-valued mapping and (X, Y, F) is a generalized approximation space. Define $F : X \to P(Y)$ $as \sim F := \{(x,y) \in X \times Y \mid y \in F(x)\}$ and for any subset A of Y , the generalized lower and upper approximations, F-(A) and F+(A), are defined by $F_-(A) = \{x \in X \mid F(x) \subseteq A\}$ and $F_+(A) = \{x \in X \mid F(x) \cap A \neq \emptyset\}$. We say that the pair F-(A),F+(A) is a generalized rough set.

Definition 3.5. Suppose that F : X →P(Y) is A set-valued mapping. We called F is a set-valued BO/ BH/Z-morphism if it satisfies :$(\forall x, y \in X) (F(x) * F(y) \subseteq F(x*y))$. A set-valued mapping t : X → P(Y) is called a strong set-valued BO/BH/Z Imorphism if it satisfies: $(\forall x, y \in X) (F(x) * F(y) = F(x * y))$.

## 4. Conclusion

This paper presents the new concepts of rough BO/BH/Z- Algebra as extended of the concept of BO/BH/Z-algebra respectively. The concept of a (strong) set-valued BO/BH/Z-morphism in BO/BH/Z algebras is investigated with several properties by Using the concept of generalized approximation space and ideal of BO/BH/Z-algebra, some properties are studied. We are sure that the results have some applications, so let us open the door to further finding new results in future work.

Acknowledgments: Acknowledgments to any discussion and suggestions from the staff of the Mathematics Department of Ajdabiya University for help in this paper. However, any comments and suggestions are acknowledged for all anonymous reviewers.